\documentclass[10pt]{article}
\usepackage{amssymb}
\usepackage{amsmath}
\usepackage{amsthm}
\usepackage{float}
\usepackage{layout}
\usepackage{amsmath}
\usepackage{graphicx}
\usepackage{graphics}
\usepackage[all]{xy}
\theoremstyle{plain}
\newtheorem{Teo}{Theorem}[section]
\newtheorem{Def}[Teo]{Definition}

\newtheorem{Rem}[Teo]{Remark}
\newcommand{\Natural}{\mathbb{N}}
\newcommand{\Integer}{\mathbb{Z}}

\begin{document}
\title{Fiboquadratic sequences and extensions of the Cassini identity raised from the study of rithmomachia}
\author{Tom\'as Guardia\footnote{Department of Mathematics. Gonzaga University. E-mail: {\tt guardia@gonzaga.edu}}
 \and Douglas Jim\'enez\footnote{Mathematics Section, Universidad Nacional Experimental Polit\'ecnica ``Antonio Jos\'e de Sucre'',  Vice-rectory of Barquisimeto. {\tt dougjim@gmail.com}}}
\maketitle
\begin{center}
\small{  \emph{To David Eugene Smith, in memoriam.}}
\end{center}
\begin{flushleft}
\footnotesize{Mathematics Subject Classification: 01A20, 01A35, 11B39 and 97A20.\\
Keywords: pythagoreanism, golden ratio, Boethius, Nicomachus, De Arithmetica, fiboquadratic sequences, Cassini's identity and rithmomachia.}
\end{flushleft}
\begin{abstract}
 In this paper, we introduce fiboquadratic sequences as a consequence of an extension to infinity of the board of rithmomachia. Fiboquadratic sequences approach the golden ratio and provide extensions of Cassini's Identity.

 \end{abstract}

\section{Introduction}
Pythagoreanism was a philosophical tradition, that left a deep influence over the Greek mathematical thought. Its path can be traced until the Middle Ages, and even to present. Among medieval scholars, which expanded the practice of the pythagoreanism, we find Anicius Manlius Severinus Boethius (480-524 A.D.)  whom by a free translation of \emph{De Institutione Arithmetica} by Nicomachus of Gerasa, preserved the pythagorean teaching inside the first universities. In fact, Boethius' book became the guide of study for excellence during quadrivium teaching, almost for 1000 years. The learning of arithmetic during the quadrivium, made necessary the practice of calculation and handling of basic mathematical operations. Surely, with the mixing of leisure with this exercise, Boethius' followers thought up a strategy game in which, besides the training of  mind calculation, it was used to preserve pythagorean traditions coming from the Greeks and medieval philosophers. Maybe this was the origin of \emph{the philosophers' game} or \emph{rithmomachia}. Rithmomachia (RM, henceforward) became the number game by excellence that picked up all pythagorean spirit; through its major victories, players can be trained in calculation of arithmetic, geometric and harmonic proportions, and its frequent practice would be concomitant to quadrivium study (arithmetic, geometry, music and astronomy). Etymologically rithmomachia is a word composed by \emph{arithmos} ($\alpha\rho\iota\vartheta\mu o\varsigma$, \emph{number}; maybe also \emph{rithmos, rhythm} as an apocope), and \emph{machia} ($\mu\alpha\chi\eta$, \emph{battle}); this leads to the meaning \emph{harmonic battle of numbers}, which perfectly synchronizes all the pythagorean essence of proportions theory.

Our main aim in this paper, is to prove that rithmomachia pieces are the first six terms of eight sequences, that we will call \emph{fiboquadratic sequences}, built up with the terms of extended Fibonacci sequence, and they hold properties of Fibonacci sequence related with the golden ratio. Also, we will prove that fiboquadratic sequences lead to extensions of Cassini's identity. We would like to think that pythagorean tradition of RM makers, arranged  on the board these fiboquadratic sequences on both numbers armies with the purpose of provide a harmonic and golden battle even though properties of a fiboquadratic sequences, should not be known by the game makers.

In another work \cite{rules} Guardia, Jim\'enez and Gonz\'alez present a set of new rules for RM, which have been tested and they have reached stability and show a complete explanation of the pieces movement, the different kinds of captures, illustrated with examples. Additionally the authors corrected the progression tables that have come to our days from several sources (medieval and modern). This work can be a starting point if someday rithmomachia get birth again and the purpose could give the first steps towards a standardization.

We have structured this work as follows: first, we describe briefly the main features of RM, starting from the pieces, how they move, the several kinds of captures and the two styles of victories. Second, we show how according to Boethian number classification, the pieces of rithmomachia are disposed on the board such that, they are the first terms of a family of fiboquadratic sequences. We will see how fiboquadratic sequences constitute a general form of some sequences, which have been studied by several authors along the last ten years. And finally, in the third part of the paper, we will show how the fiboquadratic sequences lead to natural extensions of Cassini's identity. We finish the work with a possible new proof of one well known identity.

\section{Rithmomachia}
Rithmomachia is a strategy board game created during the 11th century, with the purpose of supporting the study of the free translation made by Boethius of \emph{De Institutione Arithmetica} by Nicomachus of Gerasa, during the quadrivium course in the first universities and monastic schools. In \cite{smith}, Smith and Eaton quoted that the game was popularized in some social elite and had no attraction to the common people. We also see in this work how the game was based on Greek number theory and the authors refer some manuscripts from eleventh, twelfth and thirteenth centuries. Some scholars attributed the game to Pythagoras, but we have no records tracing of RM until the Ancient Ages. So, doubtless RM is a board strategy medieval game (See, Moyer \cite{ann}, page 2).

The pieces are rounds, triangles and squares. There is an special piece, a pyramid, made of pieces of above shapes, in decreasing sizes. The pieces shapes agree with pythagoreans conceptions about numbers. With the exception of pyramids, each piece has two faces with the same number and different color.  As we shall see in the next sections, the numerical distribution of the pieces seems to be arbitrary, but, they are really the first terms of a family of sequences, that in this paper, we introduce as \emph{fiboquadratic sequences}.

For a detailed explanations of the rules, the moving of the pieces, illustrated with examples and complete descriptions of minors and majors victories, we refer the reader to \cite{rules}. As we have said, the aim of this work is to show how fiboquadratic sequences naturally arise in the context of RM, as we shall show in the next sections. On figure \ref{Boissiere}, we have the classical diagram of de Boissi\`ere (coming up from Renaissance) that shows the initial positions of the pieces at the game start. The Venezuelan Club of Rithmomachia, an extension group of the Faculty of Science of the Central University of Venezuela, has been working in the promotion, diffusion and recovering of the healthy practice of RM, with the certainty that RM is a powerful tool to improve the mathematical learning of small calculations. The Club have tested the rules presented in \cite{rules} in the last two years and has verified their stability. Today, the first steps have been giving to set up the Gonzaga Rithmomachia Club, as an Academic Club of the Gonzaga University, which will have the main goal to spread the practice and the academic-scientific study of rithmomachia in the United States of America.
\begin{center}
\includegraphics[width=4cm,height=8cm]{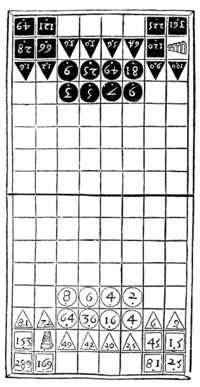}
\begin{figure*}[h]
\caption{\label{Boissiere}\small  Claude de Boissi\`ere's diagram from 1556.}
\end{figure*}
\end{center}

\section{Fiboquadratic sequences on rithmomachia}
\subsection{Nicomachus and Boethius}
Anicius Manlius Severinus Bo\"ethius (or Boethius) (480-524 A.D.) is considered the last philosopher of the antiquity  and the first philosopher of the Middle Ages. As pythagorean, he is well known by his work \emph{De consolatione philosophiae}, composed  during the last years of his life, being in prison. It is a treatise written in verse and prose, the book is a confession of his anguish by the captivity. Another work of Boethius --more important to our subject-- was a free translation of  \emph{De Institutione Arithmetica} by Nicomachus of Gerasa. In this book, we find the current classification of numbers into even and odd numbers, according to  the form of their factors and according to the sum of their divisors, each definition admits an opposite one.

\subsection{Rithmomachia pieces according to Boethius}
Mathematically speaking,  an interesting fact of this research, is  that we find inside \emph{De Institutione Arithmetica} inequalities between integers numbers that RM creators used for the numerical distribution of the pieces. In the next definition we give such number classifications (For more details, we refer the reader to \cite{boecio}). Let $(F_{n})_{n\geq 0}$ the Fibonacci numbers, that is: $F_0=0$, $F_1=1$ and $F_n=F_{n-1}+F_{n-2}$ if $n\geq 2$. The extended Fibonacci numbers are $F_{-n}=(-1)^{n+1}F_n$ if $n\geq 1$. 

\begin{Def}[Nicomachus and Boethius]\label{numeros}\rm
  A number $n$ (positive integer) is a \emph{multiple} of $m$, if $n=km$ for some $k\in\Natural$. A number $p$ (positive integer) is a \emph{superparticular} of $m$  if  $p=m(1+\frac{1}{n})$, where $n$ is a proper divisor of $m$. We say that $p$ is a \emph{superpartient} of $m$ if $p=m(1+\frac{k}{n})$,  where $n$ is a proper divisor of $m$ and $k=2,\cdots,n-1$. A number $p$ is \emph{multiple superparticular} of $m$ if $p=m(r+\frac{1}{n})$, where $n$ is a proper divisor of $m$ and $r > 1$ and it is \emph{multiple superpartient} if $p=m(r+\frac{k}{n})$ where $n$ is a proper divisor of $m$ and $k=2,\cdots,n-1$.
\end{Def}

\begin{Rem}\rm
 These definitions were given by Boethius in the context of naturals numbers, but we can extend them to integers numbers, as usual. The reader may find a more detailed explanation in Heath \cite{heath}, chapter III, pages 101-104.
 \end{Rem}

 The numerical distribution of the pieces in RM follows the Boethian classification of numbers. Table  \ref{rmboeciana} shows the values of the pieces of the white army on the left and the black army on the right.

\begin{table*}[h!]
\renewcommand{\arraystretch}{1.3}
\[\begin{tabular}{|l|cccc|cccc|}
\multicolumn{1}{c}{\ } & \multicolumn{4}{c}{\text{White army}} & \multicolumn{4}{c}{\text{Black army}}\\
\hline

$p=n$ & $2$ &   $4$ &   $6$ &   $8$      & $3$ & $5$ & $7$ & $9$   \\

$q=np$&$4$ & $16$ & $36$ & $64$      & $9$ & $25$ & $49$ & $81$ \\

$r=(1+\frac{1}{n})q$&$6$ & $20$ & $42$ & $72$      & $12$ & $30$ & $56$ & $90$ \\

$s=(1+\frac{1}{n})r$& $9$ & $25$ & $49$ & $81$     & $16$ & $36$ & $64$ & $100$ \\

$t=(1+\frac{n}{n+1})s$&  $15$ & $45$ & $91$ & $153$     & $28$ & $66$ & $120$ & $190$ \\

$u=(1+\frac{n}{n+1})t$&  $25$ & $81$ & $169$ & $289$    & $49$ & $121$ & $225$ & $361$\\
\hline
\end{tabular}\]
\caption{\small\label{rmboeciana} Numerical distribution of RM according to definition (\ref{numeros}).}
\end{table*}

The first row consists of the four first even and odd numbers respectively. The second row is a multiple of the first row, in fact it is the square. The third row is a  superparticular of the form $(1+\frac{1}{n})$ of the second row. The fourth row is a superparticular $(1+\frac{1}{n})$ of the third row. The fifth row is a superpartient of the form $(1+\frac{n}{n+1})$ of the fourth row and the sixth is a superpartient $(1+\frac{n}{n+1})$ of the fifth row.

A simple calculation shows that if $p=n$ then $q=np=n^2$, $r=(1+\frac{1}{n})q=(1+\frac{1}{n})n^2=n(n+1)$, $s=(1+\frac{1}{n})r=(1+\frac{1}{n})n(n+1)=(n+1)^2$, $t
=(1+\frac{n}{n+1})s=(1+\frac{n}{n+1})(n+1)^2=(n+1)(2n+1)$ and $u=(1+\frac{n}{n+1})t=(1+\frac{n}{n+1})(n+1)(2n+1)=(2n+1)^2$. So, table \ref{rmboeciana} corresponds to the polynomial pattern shown in table \ref{rmalgebraica}.
\begin{table*}[h!]
\[\renewcommand{\arraystretch}{1.3}
\begin{tabular}{|l|cccc|cccc|}
\multicolumn{1}{c}{\ } & \multicolumn{4}{c}{\text{White army}} & \multicolumn{4}{c}{\text{Black army}}\\
\hline
$p=n$&$2$ &   $4$ &   $6$ &   $8$      & $3$ & $5$ & $7$ & $9$   \\

$q=n^2$&$4$ & $16$ & $36$ & $64$      & $9$ & $25$ & $49$ & $81$ \\

$r=n(n+1)$&$6$ & $20$ & $42$ & $72$      & $12$ & $30$ & $56$ & $90$ \\

$s=(n+1)^2$&$9$ & $25$ & $49$ & $81$      & $16$ & $36$ & $64$ & $100$ \\

$t=(n+1)(2n+1)$&$15$ & $45$ & $91$ & $153$     & $28$ & $66$ & $120$ & $190$ \\

$u=(2n+1)^2$&  $25$ & $81$ & $169$ & $289$      & $49$ & $121$ & $225$ & $361$\\
\hline
\end{tabular}
\]
\caption{\small\label{rmalgebraica} Polynomial distribution of RM.}
\end{table*}

\subsection{Rithmomachia number extension to infinity}
A careful view of  table  \ref{rmalgebraica} shows an algebraic arrangement that distributes pairs of rows, such that the first row is of the form $k(n)l(n)$ and the second  is of the form $\left(l(n)\right)^2$, where  $k(n)$ and $l(n)$ are affine maps with integers coefficients. (In the first line: $k(n)=0\cdot n+1$ and $l(n)=1\cdot n+0$.) As we stated before, in the construction of table \ref{rmalgebraica} the second row is a multiple of the first row of the form $k(n)l(n)n$. The third row is a superparticular of the second row of the form $l(n)^2(1+\frac{1}{n})$, and the fourth row is a superparticular of the third row of the form $k(n)l(n)(1+\frac{1}{n})$. The fifth row is a superpartient of the fourth row of the form $l(n)^2(1+\frac{n}{1+n})$ and the sixth row is a superpartient of the fifth row of the form $k(n)l(n)(1+\frac{n}{1+n})$. If we want to continue this procedure, it is necessary to define a sequence that inductively generates pairs of rows of an extension to infinity of table \ref{rmalgebraica}. Now, how could we choose the coefficients of $k$ and $l$, if we would like to extend to infinity the rows of table \ref{rmalgebraica} keeping this same initial disposition? An incidental view helps to decide: \emph{the coefficients of the terms on $n$ that we have until now are the first terms of Fibonacci sequence}: $0$, $1$, $1$, $2$. This motivate us to define $(b_{mn})_{m\geq 0}$ as the next sequence:
\begin{equation}\label{sucesion bmn}
  b_{mn}=1+\frac{F_{m-1}n+F_{m-2}}{F_mn+F_{m-1}}
 \end{equation}
where $F_m$ (for $m\geq 0$) are Fibonacci numbers. Then, we define an infinite extension of table \ref{rmalgebraica} as:
\begin{equation}\label{sucesion amn}
  a_{0n}=1,\quad a_{mn}=\begin{cases}
    a_{m-1,n}b_{\frac{m-1}{2},n}, &\text{ if $m$ is odd }\\
    a_{m-1,n}b_{\frac{m}{2}-1,n}, &\text{ if $m$ is even }\\
  \end{cases}
\end{equation}

We have to check that the first terms of (\ref{sucesion amn}) are exactly the rows of table \ref{rmalgebraica}. In fact, if $m\geq1$ then:
$a_{1n}=a_{0n}b_{0n}=1\cdot n=n$, $a_{2n}=a_{1n}b_{0n}=n\cdot n=n^2$, $a_{3n}=a_{2n}b_{1n}=n^2(1+\frac{1}{n})=n(n+1)$, $a_{4n}=a_{3n}b_{1n}=n(n+1)(1+\frac{1}{n})=(n+1)^2$, $a_{5n}=a_{4n}b_{2n}=(n+1)^2(1+\frac{n}{n+1})=(n+1)(2n+1)$, $a_{6n}=a_{5n}b_{2n}=(n+1)(2n+1)(1+\frac{n}{n+1})=(2n+1)^2$.

If we continue this procedure to infinity, we obtain an infinite extension of table \ref{rmalgebraica}. An explicit form of (\ref{sucesion amn}) is given in the next theorem.
\begin{Teo}\label{sucesion amn explicita}
Let $(b_{mn})_{m\geq0}$ be the sequence defined in (\ref{sucesion bmn}). Then the sequence $(a_{mn})_{m\geq0}$ defined in (\ref{sucesion amn}) is:
\begin{equation}\label{rmfibonacci}
  a_{mn}=\begin{cases}
    \left(F_{\frac{m-1}{2}}n+F_{\frac{m-1}{2}-1}\right)\left(F_{\frac{m-1}{2}+1}n+F_\frac{m-1}{2}\right), &\text{if $m$ is odd} \\
    \left(F_{\frac{m}{2}}n+F_{\frac{m}{2}-1}\right)^2, &\text{if $m$ is even}
  \end{cases}
\end{equation}
where $F_k$ for each $k\in\Integer$ is an extended Fibonacci number.
\end{Teo}
\emph{Proof:} If $m=0,1$, then $a_{0n}=(F_0n+F_{-1})^2=1$ and  $a_{1n}=(F_0n+F_{-1})(F_1n+F_0)=n$. So,  (\ref{rmfibonacci}) holds. Suppose that (\ref{rmfibonacci}) holds until $m$. If $m$ is odd then $m+1$ is even, Therefore, by the properties of Fibonacci numbers, we have:
\begin{align*}
  a_{(m+1),n}&=a_{mn}b_{\frac{m+1}{2}-1,n}\\
  &=a_{mn}b_{\frac{m-1}{2},n}\\
  &=\left(F_{\frac{m-1}{2}}n+F_{\frac{m-1}{2}-1}\right)\left(F_{\frac{m-1}{2}+1}n+F_{\frac{m-1}{2}}\right)\left(1+\frac{F_{\frac{m-1}{2}-1}n+F_{\frac{m-1}{2}-2}}{F_{\frac{m-1}{2}}n+F_{\frac{m-1}{2}-1}}\right)\\
  &=\left(F_{\frac{m+1}{2}}n+F_{\frac{m+1}{2}-1}\right)^2
\end{align*}
Otherwise, if $m$ is even then $m+1$ is odd and again:
\begin{align*}
  a_{(m+1),n}&=a_{mn}b_{\frac{m}{2},n}\\
  &=\left(F_{\frac{m}{2}}n+F_{\frac{m}{2}-1}\right)^2\left(1+\frac{F_{\frac{m}{2}-1}n+F_{\frac{m}{2}-2}}{F_{\frac{m}{2}}n+F_{\frac{m}{2}-1}}\right)\\
  &=\left(F_{\frac{m}{2}}n+F_{\frac{m}{2}-1}\right)\left(F_{\frac{m}{2}+1}n+F_{\frac{m}{2}}\right)
\end{align*}
By induction we have that (\ref{rmfibonacci}) holds for every $m\geq 0$ and the theorem is proved. \qed

The equation (\ref{rmfibonacci}) is an infinite extension of table \ref{rmalgebraica}, as it is shown in table \ref{rminfinita}, which it is increasingly ordered with respect to the columns.

\begin{center}
\begin{tabular}{lccccccccc}


 $a_{mn}$ & $a_{m2}$ & $a_{m3}$ & $a_{m4}$ & $a_{m5}$  &$a_{m6}$ & $a_{m7}$ & $a_{m8}$ & $a_{m9}$ \\

1& 1& 1& 1& 1& 1& 1& 1&1\\
$n$& 2& 3& 4& 5& 6& 7& 8& 9\\
$n^2$& 4& 9& 16& 25& 36& 49& 64& 81\\
$n(n + 1)$& 6& 12& 20& 30& 42& 56& 72& 90\\
$(n + 1)^2$& 9& 16& 25& 36& 49& 64& 81& 100\\
$(n + 1)(2n + 1)$&  15& 28& 45& 66& 91& 120& 153& 190\\
$(2n + 1)^2$& 25& 49& 81& 121& 169& 225& 289& 361\\
$(2n + 1)(3n + 2)$& 40& 77& 126& 187& 260& 345& 442& 551\\
$(3n + 2)^2$& 64& 121& 196& 289& 400& 529& 676& 841\\
$(3n + 2)(5n + 3)$& 104& 198& 322& 476& 660& 874& 1118& 1392\\
$(5n + 3)^2$& 169& 324& 529& 784& 1089& 1444& 1849& 2304\\
$(5n + 3)(8n + 5)$& 273& 522& 851& 1260& 1749& 2318& 2967& 3696\\
$(8n + 5)^2$& 441& 841& 1369& 2025& 2809& 3721& 4761& 5929\\
$(8n + 5)(13n + 8)$& 714& 1363& 2220& 3285 &4558 &6039 &7728 &9625\\
$(13n + 8)^2$&1156& 2209& 3600& 5329& 7396& 9801& 12544& 15625\\
$\vdots$&  $\vdots$ & $\vdots$ & $\vdots$ & $\vdots$     & $\vdots$& $\vdots$ & $\vdots$ & $\vdots$ \\
\end{tabular}
\begin{table}[h]
\caption{\label{rminfinita}\small An infinite extension of RM.}
\end{table}
\end{center}

Following rithmomachia we use for the columns the natural numbers from 2 until 9, but we can extend this feature until desired, so table \ref{rminfinita} may be extended to infinity by the columns.

Notice that on RM, we do not only get the numerical equality, but the geometrical equality. Since the superposition of two triangles is a parallelogram and since a row of squares is preceded by two rows of triangles, we observe that the creators of RM respected this spirit in adding two triangles to obtain one square.  We realized this is a recursive definition which motivated us to extend the board to infinity and find some interesting relations among numbers that we shall explain on the next section.

\begin{center}
\[
\begin{matrix}
\includegraphics[width=.5cm,height=.5cm]{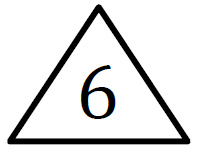}&+&
\includegraphics[width=.5cm,height=.5cm]{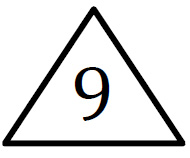}&=&
\includegraphics[width=.5cm,height=.5cm]{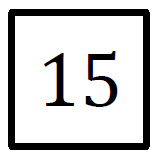} \\
\includegraphics[width=.5cm,height=.5cm]{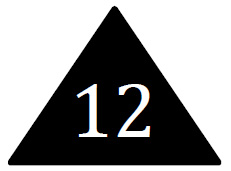}&+&
\includegraphics[width=.5cm,height=.5cm]{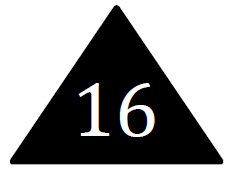}&=&
\includegraphics[width=.5cm,height=.5cm]{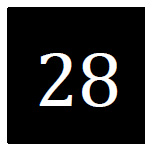} \\
\end{matrix}
\]
\end{center}
\begin{figure*}[h]
\caption{\small The geometrical sum of two triangles is an square.}
\end{figure*}

If we can imagine an extension \emph{ad infinitum} of the board we will find two armies of sequences ready to start  the \emph{harmonic battle of numbers}. Doubtlessly, \emph{the philosophers' game}, was at its time, a complete expression of the main elements of pythagoreanism, that idealized an aesthetic universe, leaded by the supreme fundamentals laws of beauty.

\subsection{Fiboquadratic sequences}

The construction of an infinite board of RM made on the previous section, resulted in a family of sequences, this lead us to introduce the following definition.

\begin{Def}\rm\label{deffibocuadratica}
For any $n\in \mathbb{N}$, the \emph{fiboquadratic sequence} generated by $n$ is the sequence $\displaystyle (a_{mn})_{m\in \mathbb{N}}$ defined as
\[
a_{mn}=\left\{
\begin{array}{lcl}
\left(F_{\frac{m-1}{2}}n+F_{\frac{m-1}{2}-1}\right)\left(F_{\frac{m-1}{2}+1}n+F_{\frac{m-1}{2}}\right), & & \text{if $m$ is odd}\\
\left(F_{\frac{m}{2}}n+F_{\frac{m}{2}-1}\right)^2, & & \text{if $m$ is even}
\end{array}
\right.
\]
where $F_n$ for each $n\in\Natural$ is a term of the extended Fibonacci sequence.
\end{Def}

\begin{Rem}\rm
This definition could be extended to the most general case where we use a real variable $t$ instead of the index $n$, that is:
\[a_m(t)=\left\{
\begin{array}{lcl}
\left(F_{\frac{m-1}{2}}t+F_{\frac{m-1}{2}-1}\right)\left(F_{\frac{m-1}{2}+1}n+F_{\frac{m-1}{2}}\right), & & \text{if $m$ is odd}\\
\left(F_{\frac{m}{2}}t+F_{\frac{m}{2}-1}\right)^2, & & \text{if $m$ is even},
\end{array}
\right.
\]
but we do not use this variation in this paper.
\end{Rem}

\begin{Rem}\rm
Every fiboquadratic sequence $(a_{mn})_{m\in\Natural}$ can be extended to negative indexes as usual; that is, switching the index from natural to integer. In this case, we say that, $(a_{kn})_{k\in\Integer}$ is the \emph{extended fiboquadratic sequence} of $(a_{mn})_{m\in\Natural}$.
\end{Rem}

\begin{Rem}\rm
If we let $n=1,3$, then the sequences $(a_{m1})_{m\in\Natural}$ and $(a_{m3})_{m\in\Natural}$ on definition (\ref{deffibocuadratica}) are respectively the sequences A006498 and  A006499 recorded in \emph{The On--Line Encyclopedia of Integer Sequences}  (See for instance, \cite{6498} and \cite{6499}). If $n=2$, the sequence $(a_{m2})_{m\in\Natural}$ is the truncated version of $(a_{m1})_{m\in\Natural}$ (see, \cite{6498}). The reader can check that the remaining six sequences $(a_{mn})_{m\in\Natural}$ with $n=4,5,6,7,8,9$ are not recorded in this encyclopedia. Hence, fiboquadratic sequences defined in (\ref{deffibocuadratica}) are natural generalizations of sequences A006498 and A006499.
\end{Rem}

As fiboquadratic sequences are defined in terms of Fibonacci numbers, we hope these sequences are related to the golden number $\alpha$ and $\beta=-1/\alpha$. We formally address this fact in the next theorem.

\begin{Teo}\label{fibocuadratica}
 The sequence of quotients of successive terms of any fiboquadratic sequence converges to the golden ratio, i. e.
 \begin{equation}
 \label{limitefibocuadratica} \lim_{m \rightarrow \infty} \frac{a_{m+1,n}}{a_{mn}}=\alpha, \qquad \forall n \in \Natural.
 \end{equation}
\end{Teo}
\emph{Proof:} A carefully examination of table 5 and definition 4.4 show us that after a change of variable, the rows of the matrix follow the ternary scheme
\begin{align*}
A(m)&=(F_mn+F_{m-1})^2\\
B(m)&=(F_mn+F_{m-1})(F_{m+1}n+F_m)\\
C(m)&=(F_{m+1}n+F_m)^2.
\end{align*}

Then
\[\frac{B(m)}{A(m)}=\frac{C(m)}{B(m)}=\frac{F_{m+1}n+F_m}{F_mn+F_{m-1}}=\frac{\frac{F_{m+1}}{F_m}n+1}{n+\frac{F_{m-1}}{F_m}}\]
but we know that $\frac{F_{m+1}}{F_m} \rightarrow \alpha$ as $m\rightarrow \infty$, therefore
\[
\lim_{m \rightarrow \infty}\frac{B(m)}{A(m)}=\lim_{m \rightarrow \infty}\frac{C(m)}{B(m)}=\lim_{m \rightarrow \infty}\frac{\frac{F_{m+1}}{F_m}n+1}{n+\frac{F_{m-1}}{F_m}}=\frac{\alpha n+1}{n+\frac{1}{\alpha}}=\alpha.
\] \qed

\begin{Rem}\rm
It is not difficult to show that
\[\lim_{m \rightarrow -\infty} \frac{a_{m+1,n}}{a_{mn}}=\beta, \qquad \forall n \in \Natural.\]
\end{Rem}

The matrix shown in table \ref{rminfinita} is actually a family of eight fiboquadratic sequences ordered by columns. But these facts can be seen in a more general view, introducing the generalized Fibonacci sequences that can be expressed in an even more general view, as in the next definition.

\begin{Def}\label{generalFibonacci}
If $a,\, b\in \Integer$,--not both of them equal $0$-- the sequence
\[G_1=a,\qquad G_2=b,\qquad G_{m+1}=G_{m} + G_{m-1},\quad m \geq 2\]
is called a generalized Fibonacci sequence.
\end{Def}

In the context of generalized Fibonacci sequence there is a very important number, the so called \emph{characteristic} of the sequence, that is:
\begin{equation}\label{characteristic}
\mu = \det\begin{pmatrix}G_3 & G_2\\ G_2 & G_1 \end{pmatrix}=a^2+ab-b^2.
\end{equation}

Then we can put our eyes in the factors of the fiboquadratic sequences because
\begin{equation}\label{def_fibogeneral_a}
G_m=F_{m-1} n + F_{m-2}, \quad m \geq 1
\end{equation}
defines a particular case of generalized Fibonacci sequence. This can be expressed as a theorem.

\begin{Teo}\label{generalFibonacciRithmo}
If $(G_m)$ is defined as in equation (\ref{def_fibogeneral_a}) then
\[G_1=1,\qquad G_2=n,\qquad G_{m+1}=G_{m} + G_{m-1},\quad m \geq 2\]
\end{Teo}

\emph{Proof:} $G_1=F_0 n + F_{-1}=1$, $G_2=F_1 n+ F_0=n$ and
\begin{align*}
G_{m} + G_{m-1} &= (F_{m-1} n + F_{m-2})+(F_{m-2} n + F_{m-3})\\
    &= (F_{m-1} + F_{m-2})n+(F_{m-2}  + F_{m-3})\\
    &= F_m n + F_{m-1}\\
    &= G_{m+1},
\end{align*}
as we expected.\qed

\section{Cassini's identities}
\subsection{Introduction}
In Dunlap \cite{dunlap},  Grimaldi \cite{grimaldi}, Koshy \cite{koshy} and Vajda \cite{vajda} we can see that there are a lot of properties of Fibonacci numbers. One of them --that will be very important for us-- was discovered in 1860 by the Italian-born French astronomer and mathematician Giovanni Domenico Cassini (1625-1712) and this result was also discovered independently in 1753 by the Scottish mathematician and landscape artist Robert Simson (1687-1768); it can be formulated as
\begin{equation}\label{cassini}
F_{m-1}F_{m+1}-F_m^2=(-1)^m,\qquad \forall m\in\Natural.
\end{equation}

Cassini's identity involves any three consecutive members of the Fibonacci sequence, but it can be extended to strings of any size of consecutive terms of the sequence; there is a relation between the product of the extreme terms of the string and the middle terms, but this relation depends upon the parity of the size of the string. If $n,\, k \in \Natural$, such that $0< k < n$, the following equations explain this idea:

\textbf{Odd size strings (Catalan's theorem)}
\begin{equation}\label{cassini_ext_impar}
F_{m+k}F_{m-k}-F_m^2 = (-1)^{m+k+1}F_k^2
\end{equation}

(Cassini's identity is the particular case $k=1$.)

\textbf{Even size strings (Vajda's Formula)}
\begin{equation}\label{cassini_ext_par}
F_{m+k+1}F_{m-k}-F_mF_{m+1} = (-1)^{m+k+1}F_kF_{k+1}
\end{equation}

A proof of equation (\ref{cassini_ext_impar}) can be found in Koshy (\cite{koshy}, page 83). In addition, equation (\ref{cassini_ext_par}) is found in Vajda (\cite{vajda}, page 28) and it can be proved similarly. If we have a generalized Fibonacci sequence as that of definition \ref{generalFibonacci} the above equations  take these forms:

\textbf{Odd size strings (Tagiuri, 1901)}
\begin{equation}\label{general_cassini_ext_impar}
G_{m+k}G_{m-k}-G_m^2 = (-1)^{m+k+1}\mu F_k^2
\end{equation}

\textbf{Even size strings (Tagiuri, 1901)}
\begin{equation}\label{general_cassini_ext_par}
G_{m+k+1}G_{m-k}-G_mG_{m+1} = (-1)^{m+k+1}\mu F_kF_{k+1}
\end{equation}

 In the next section, we will show that fiboquadratic sequences bring within them a natural extension of the Cassini's Identity.

\subsection{Cassini's identities on fiboquadratic sequences}\label{cassinirm}
If one look forward to  table \ref{rminfinita}, it should be noticed that the sum of any two consecutive rows follow a double pattern: (a) an even row added to the next odd row gives the following even row, (b) odd row with the next even row is, up to a fixed constant that depends on the column, the next odd row. Let us see this last case, column by column, showing with the first three of them:
\[
\begin{array}{rrrc}
1+2=4-\boldsymbol{1},&4+6=9+\boldsymbol{1},&9+15=25-\boldsymbol{1},&\ldots\\
1+3=9-\boldsymbol{5}, & 9+12=16+\boldsymbol{5}, & 16+28=49-\boldsymbol{5}, & \ldots\\
1+4=16-\boldsymbol{11}, & 16+20=25+\boldsymbol{11}, & 25+45=81-\boldsymbol{11}, & \ldots\\
1+5=25-\boldsymbol{19}, & 25+30=36+\boldsymbol{19}, & 36+66=121-\boldsymbol{19}, & \ldots
\end{array}
\]

The next table shows the fixed constants for the eight fiboquadratic sequences founded on the infinite extension of RM
\begin{center}
\begin{tabular}{|c|c|}\hline
    \text{Column Number} &\text{Fixed Constant}\\\hline
    2&1\\\hline
    3&5\\\hline
    4&11\\\hline
    5&19\\\hline
    6&29\\\hline
    7&41\\\hline
    8&55\\\hline
    9&71\\\hline
\end{tabular}
\begin{table}[h]
\caption{\label{constantesfijas}\small The fixed constants arising from the sum of two next rows of table \ref{rminfinita}.}
\end{table}
\end{center}

The reader can check that on the others columns, the sum of the corresponding rows differs up the fixed constants shown on the table \ref{constantesfijas}.

As an example, the behavior described above over the first column of table \ref{rminfinita} leads us to conjecture the next equation.
\begin{equation}\label{conjetura}
(2F_m+F_{m-1})^2+(2F_m+F_{m-1})(2F_{m+1}+F_m)=(2F_{m+1}+F_{m})^2-(-1)^m.
\end{equation}

The definition of the Fibonacci sequence allows us to rewrite (\ref{conjetura}) as
\begin{equation}\label{cassiniforma1}
  F_{m+2}F_{m+4}=F_{m+3}^2-(-1)^m.
\end{equation}
After a change of variable, we note that, in fact (\ref{cassiniforma1}) is Cassini's identity, as show in (\ref{cassini})
\[
F_{m-1}F_{m+1}-F_m^2=(-1)^{m}
\]
So, reversing the steps, it is clear that equation (\ref{conjetura}) holds.

A closed examination of table \ref{constantesfijas} let us to conjecture
\[c_0=-1, \qquad c_n=c_{n-1}+2(n-1),\quad n \geq 1\]
as a formula for the fixed constants of the table; this is a non--homogeneous recurrence relation equivalent to
\begin{equation}\label{polifibonacci}
c_n = n^2-n -1, \quad n \geq 1
\end{equation}
as can be seen easily by induction. The right side of equation (\ref{polifibonacci}) is the well known \emph{Fibonacci's polynomial}. But
\[c_n = -\det \begin{pmatrix}
n+1 & n \\ n & 1
\end{pmatrix}=-\mu\]
is the negative of the characteristic of the generalized Fibonacci sequences defined in equation (\ref{def_fibogeneral_a}). These heuristic reasonings let us to think that if $(a_{mn})$ is a fiboquadratic sequence then
\begin{equation}\label{teoremaestrella}
\begin{array}{lcl}
a_{mn}+a_{m+1,n}=a_{m+2,n}, & & \text{if $m$ is even}\\
a_{mn}+a_{m+1,n}=a_{m+2,n}+(-1)^{\frac{m+1}{2}}(n^2-n-1), & & \text{if $m$ is odd}
\end{array}
\end{equation}
but again we can put the problem in a more general view.

For this purpose we can define a general fiboquadratic sequence upon any generalized Fibonacci sequence $(G_m)$, in the following way:
\begin{equation}\label{generalFiboquad}
(a_m)=(G_1^2,\, G_1G_2,\, G_2^2,\, G_2G_3,\, G_3^2,\, \ldots),
\end{equation}
so we can state our general and central theorem:

\begin{Teo}\label{teoremageneral}
If $(a_{m})$ is any fiboquadratic sequence then
\[
\begin{array}{lcl}
a_{m}+a_{m+1}=a_{m+2}, & & \text{if $m$ is even}\\
a_{m}+a_{m+1}=a_{m+2}+(-1)^{\frac{m+1}{2}}\mu, & & \text{if $m$ is odd}
\end{array}
\]
\end{Teo}

\emph{Proof:} In the definition of $(a_{m})$ we can define two important subsequences:
\[(o_m)=(a_{2m-1})=(G_1^2,\, G_2^2,\,  G_3^2,\, \ldots)\]
of the elements of $(a_{m})$ with odd indexes, and
\[(e_m)=(a_{2m})=(G_1G_2,\,  G_2G_3,\, G_3G_4,\, \ldots)\]
of the elements of $(a_{m})$ with even indexes.

Then, the left side of the first equation can be rewritten as
\begin{align*}
e_{m-1} + o_m &= G_{m-1}G_{m}+G_{m}^2\\
   &= G_{m}(G_{m-1}+G_{m})\\
   &= G_{m}G_{m+1}\\
   &= e_{m},
\end{align*}
and the left side of the second, as
\begin{align*}
o_m + e_m &= G_{m}^2+G_mG_{m+1}\\
   &= G_{m}(G_m+G_{m+1})\\
   &= G_{m}G_{m+2}\\
   &= G_{m+1}^2+(-1)^{m+1}\mu\tag{By Eq. (\ref{general_cassini_ext_impar})}\\
   &= o_{m+1}+(-1)^{m+1}\mu.
\end{align*}
but both results are the contents of the theorem.\qed

So, equations (\ref{teoremaestrella}) are the particular case of this theorem for the fiboquadratic sequences defined over generalized Fibonacci sequences of theorem \ref{generalFibonacciRithmo}. We want remark, that besides the convergence to the golden ratio, the general Cassini's identity naturally appear in the context of fiboquadratic sequences. In fact,  after a change of variable, the general Cassini's identity $G_{m+1}G_{m-1}=G_m^2+(-1)^m(n^2-n-1)$ is the second equation in (\ref{teoremaestrella}).

\subsection{One last formula}

The following formula -concerning the usual Fibonacci sequence- is well known:
\begin{equation}\label{sumafibocuadrados}
F_1^2+ F_2^2+\cdots + F_n^2= F_nF_{n+1},
\end{equation}
but we want to show a little variation of this, derived from the considerations over theorem \ref{teoremageneral}.

Equation $e_{m-1}+o_m = e_{m}$ in the proof of theorem  \ref{teoremageneral} can be written as $e_m=e_{m-1}+G_m^2$ and reversing the recursive property, we can reach
\begin{equation}\label{retroceso}
e_m = e_1 + G_2^2 + \cdots + G_m^2
\end{equation}

Substituting equation (\ref{retroceso}) in $o_m + e_m = o_{m+1}+(-1)^{m+1}\mu$ from the same theorem and after application of definitions we arrive to:
\[G_m^2 + G_1G_2 + G_2^2 + \cdots + G_m^2 = G_{m+1}^2 + (-1)^{m+1}\mu,\]
which --by $G_1=a$, $G_2=b$-- can be rewritten as
\[ab + G_2^2 + \cdots + G_m^2 = G_{m+1}^2 - G_m^2 + (-1)^{m+1}\mu,\]
or
\[(ab-a^2) + G_1^2 + G_2^2 + \cdots + G_m^2 = (G_{m+1} - G_m)(G_{m+1} + G_m) + (-1)^{m+1}\mu,\]
therefore
\[G_1^2 + G_2^2 + \cdots + G_m^2 = a(a-b) + G_{m-1}G_{m+2} + (-1)^{m+1}\mu,\]
that --after application of (\ref{general_cassini_ext_par})-- becomes
\[G_1^2 + G_2^2 + \cdots + G_m^2 = a(a-b) + G_{m}G_{m+1} + (-1)^m\mu + (-1)^{m+1}\mu,\]
and, finally,
\begin{equation}\label{sumafibocuadradosgeneral}
G_1^2 + G_2^2 + \cdots + G_m^2 = a(a-b) + G_{m}G_{m+1},
\end{equation}
a generalization of equation (\ref{sumafibocuadrados}) proved in a possible new way derived from the study of a medieval game.

\section{Conclusions}
The pythagoreans conceived the universe within an order and harmony ruled by numbers. We have shown that we can find this order and harmony inside rithmomachia. Rithmomachia was used during the Middle Ages for the exercise of arithmetics, geometry and music. In fact, the excellentissima victory given by the quartet $(4,6,8,12)$ has into their irreducible fractions, harmonics of the pentatonic pythagorean music scale (see \cite{smith}, page 79.). So, in order to rescue the healthy practice of RM  as a  discipline for the learning could be interesting as a pedagogical aim. Since 2013, the \emph{Venezuelan Rithmomachia Club} has dedicated to promote the study and practice of RM with highly grateful and satisfactory results. Now, is in progress the \emph{Gonzaga Rithmomachia Club} and this Academic Organization will be main promoter and diffuser of rithmomachia in the United States of America. With this research we have evidenced that rithmomachia makers fully understood the pythagoreans principles of numbers and gave --probably unknowing it-- the first steps in finding an approximation of the golden ratio.

We have introduced fiboquadratic sequences, as an heuristic construction of an infinite extension of the board of RM. As a consequence, we have shown these fiboquadratic sequences have a close behavior to the Fibonacci sequence --because, successive quotients of a fiboquadratic sequence approach to the golden ratio and they provide a generalization of Cassini's identity-- also, fiboquadratic sequences are generalizations of sequences A006498 and A006499 of the On Line Encyclopedia of Integer Sequences. In fact, they are sequences $(a_{m2})_{m\in\Natural}$ and $(a_{m3})_{m\in\Natural}$ defined in table \ref{rminfinita}. Several authors have found the generating functions, the recurrence relation and some combinatorial interpretations of both sequences. For example, A006498 counts the number of compositions of $n$ with 1's 3's and 4's and A006499 is the number of restricted circular combinations (for more details, we refer the reader to \cite{6498} and \cite{6499}). We think that with this work we open a research line to find new properties of a general fiboquadratic sequence.

\section{Acknowledgements}
Both authors thank the encouragement and noteworthy comments of Ann Elizabeth Moyer and Sophie Caflisch which helped to improve the paper. And to the master David Eugene Smith for his visionary role, who inspired us to explore inside the unknown and discover beautiful properties of numbers.

\end{document}